\documentclass[11pt,oneside]{amsart}

\usepackage{amsmath,graphicx,amssymb,latexsym}

\def\R{\mathbb{R}}
\def\N{\mathbb{N}}
\def\Z{\mathbb{Z}}
\def\G{G_3}
\def\B{B_3}
\def\Geods{\Omega}

\newtheorem{thm}{Theorem}[section]
\newtheorem{lemma}[thm]{Lemma}
\newtheorem{prop}[thm]{Proposition}
\newtheorem{corollary}[thm]{Corollary}
\theoremstyle{definition}
\newtheorem{defn}[thm]{Definition}

\begin{document}

\title[Geodesics in $\B$]{Geodesics in the Braid Group on Three Strands}
\author[L.Sabalka]{Lucas~Sabalka$\!\,^1$}
      \address{Deptartment of Mathematics\\
               University of Illinois at Urbana-Champaign \\
               Champaign, IL 61820}
      \email{sabalka@math.uiuc.edu}

\footnotetext[1]{This work was done under the supervision of Drs.
Susan Hermiller and John Meakin at the University of
Nebraska-Lincoln, with funding from National Science Foundation
grant DMS-0071037 and from a Pepsi UCARE grant.}

\begin{abstract}
We study the geodesic growth series of the braid group on three
strands, $\B := \langle a,b|aba = bab \rangle$. We show that the
set of geodesics of $\B$ with respect to the generating set $S :=
\{a,b\}^{\pm 1}$ is a regular language, and we provide an explicit
computation of the geodesic growth series with respect to this set
of generators.  In the process, we give a necessary and sufficient
condition for a freely reduced word $w \in S^*$ to be geodesic in
$\B$ with respect to $S$. Also, we show that the translation
length with respect to $S$ of any element in $\B$ is an integer.
\end{abstract}

\maketitle

\section{Introduction}

For a finitely generated group $G$ the Geodesic Growth Series of
$G$ with respect to a generating set $S=S^{-1}$ is defined by

$$\mathcal{G}_{(G,S)}(x) = \sum_{n=0}^\infty a_nx^n,$$
where $a_n$ is the number $d_S$-geodesic words of length $n$ in
$S^*$ and where $d_S$ is the word metric on $G$ defined by $S$.
Similarly, the Spherical Growth Series is defined by

$$\mathcal{S}_{(G,S)}(x) = \sum_{n=0}^\infty b_nx^n,$$
where $b_n$ is defined to be the number of elements of $G$ at
$d_S$-distance $n$ from $1$. A power series is \emph{rational} if
it may be expressed as a quotient of two integer-coefficient
polynomials in the ring of formal power series $\Z[[x]]$. It is
well known that the growth series of a regular language is
rational \cite{ECHLPT}.

Note that regularity of the set of geodesics depends on the
generating set used (\cite{NeumannShapiro}, attributed to Cannon).
It is also known that the rationality of the spherical growth
series depends on the generating set used \cite{Stoll}.

Which groups have rational \emph{geodesic} growth series? There
are many groups which are known to have rational \emph{spherical}
growth series. Groups which have rational geodesic growth series
are less well understood, but there are still many examples. It is
known, for instance, that all discrete co-compact hyperbolic
groups \cite{Cannon} have rational geodesic growth series. More
generally, Neumann and Shapiro proved all geometrically finite
hyperbolic groups have rational geodesic growth series
\cite{NeumannShapiro}, and all word-hyperbolic groups have
rational geodesic growth series, independent of generating set
(see \cite{ECHLPT}; attributed there to Gromov and others). Also,
Loeffler, Meier, and Worthington proved this fact for right-angled
Artin and Coxeter groups \cite{LoefflerMeierWorthington}, and
Noskov showed that irreducible affine Coxeter groups have rational
geodesic growth series \cite{Noskov}. Note these references
actually show the regularity of the set of geodesics, which
implies rationality of the geodesic growth series.

In this paper, we analyze the geodesic growth series of the braid
group on three strands, $\B := \langle a,b | aba = bab \rangle$
(for a detailed description of $\B$, see \cite{ECHLPT}).  We prove
that $\B$ has a rational geodesic growth series with respect to
the generating set $S = \{a,b\}^{\pm 1}$ by showing the geodesics
are a regular language. Moreover, we provide an explicit
computation of the geodesic growth series with respect to $S$.

In fact, it has been shown (independent of this work) that $\B$
has a rational geodesic growth series, although the actual
geodesic growth series has not been computed. This was done by REU
students Griffing, Patlovany, and Talley under the guidance of Jon
McCammond \cite{McCammond}. Charney and Meier have also recently
proven that the geodesic growth series is rational for all Garside
groups (which includes the braid groups), but on different
generating sets than those considered here \cite{CharneyMeier}.

We prove the following theorems:

\begin{thm}
\label{thm:geodesics} A freely reduced word $w \in S^*$ is a
geodesic for $\B$ if and only if $w$ $w$ does not contain as
subwords any of the following:
\begin{itemize}
\item elements of both $\{ab,ba\}$ and $\{AB,BA\}$;
\item both $aba$ and either $A$ or $B$;
\item both $ABA$ and either $a$ or $b$.
\end{itemize}
\end{thm}

\begin{thm} \label{thm:GGS_B} The set of geodesics for $\B := \langle
a,b | aba = bab\rangle$ with respect to the generating set $S =
\{a,b\}^{\pm 1}$ is regular, and the geodesic growth series of
$\B$ is
$$ \mathcal{G}_{(\B,S)}(x) = \frac{x^4+3x^3+x+1}{(x^2+x-1)(x^2+2x-1)}.$$
\end{thm}

It is not known whether the geodesic growth series is rational for
all braid groups.

This paper is organized as follows.  In \S 2, background
definitions are given. In section \S 3 we introduce $\B$ and
construct its Cayley graph. In section \S 4 we provide an analysis
of possible geodesics in $\B$ to prove Theorem
\ref{thm:geodesics}, and we prove the translation length with
respect to $\{a,b,\}^{\pm 1}$ of any element of $\B$ is an
integer. In section \S 5 we give the proof of Theorem
\ref{thm:GGS_B}.

I would particularly like to thank Susan Hermiller and John Meakin
as well as Ilya Kapovich for all of their time, advice, and help.

\section{Background}

Let $G$ be a group with a generating set $S = S^{-1}$. For any
word $w \in S^*$ we denote by $\overline{w}$ the element of $G$
represented by $w$.  If $w,v \in S^*$ are two words such that
$\overline{w} = \overline{v}$, we write $w \equiv v$. Given a
nontrivial freely reduced word $w \in S^*$, $w$ may be written in
the form $w = w_1w_2...w_m$, where each $w_j$ is of the form $s^k$
for some $s\in S$ and $k \in \Z, k \neq 0$, and where for $1 \leq
j < m$, if $w_j = s^k$ then $w_{j+1} = t^{l}$ for some $l \neq 0$
with $t \neq s$. We call $w_j$ the $j$th \emph{syllable} of $w$.
For any $x \in S$, we often use capitalization, $X$, to denote
$x^{-1}$. If $w$ is a word or a path, we denote the length of $w$
by $|w|$.

Throughout this paper, we will use the term \emph{graph} to mean a
directed edge-labelled graph. Denote the set of vertices of a
graph $M$ by $VM$ and the set of edges by $EM$.  Also, let an edge
$e \in EM$ from $v \in VM$ to $w \in VM$ labelled by $s$ be
denoted $e=[v,w,s]$.

The main tool we will use to determine the geodesic growth series
for $\B$ is its Cayley graph. The \emph{Cayley graph} of a group
$G$ with respect to a generating set $S=S^{-1}$ is denoted
$C_{(G,S)}$. In this paper we will only work with the situation
when $S$ is partitioned as $S = S_0 \sqcup S_0^{-1}$.  We will
refer to elements of $S_0$ as \emph{positive labels}.  If an edge
is labelled with an element of $S_0$, it is called \emph{positive}
(or \emph{negative} otherwise).  Similarly, if a syllable of a
word is an element of $S_0$ raised to a positive power, or an
element of $S_0^{-1}$ raised to a negative power, it is also
called positive. If a syllable is raised to an even (respectively,
odd, negative, or positive) power, we say the syllable is
\emph{even} (respectively, \emph{odd}, \emph{negative},
\emph{positive}). When drawing or defining a Cayley graph we will
normally just indicate the positively labelled edges. If a
specific generating set is understood, we simply write $C_G$ for
the Cayley graph of $G$. In a Cayley graph $C = C_{(G,S)}$, a
vertex $v \in VC$ is \emph{represented by} a word $w \in S^*$ if
there is a path $p$ with label $w$ beginning at $1$ and ending at
$v$ in $C$. In this case, we do not distinguish between $p$ and
its label $w$. For any vertex $v \in VC$ represented by $w'$ in
$S^*$ and any word $w \in S^*$, we denote by $vw$ the vertex in
$VC$ represented by $w'w$.

If $p$ is of minimal length among all paths with the same end
points, then $p$ is called \emph{geodesic} in $C_G$. Similarly, a
word $w \in S^*$ which is of minimal length among all words $w'
\in S^*$ for which $w' \equiv w$ is called \emph{geodesic} in $G$
(with respect to $S$) . For a given element $g \in G$, if the word
$w \in S^*$ is such that $\overline{w} = g$ and $w$ is a geodesic
in $G$, then the \emph{cone type} of $g$, and of the vertex in the
Cayley graph which represents $g$, is defined as the set of all
strings $\gamma$ such that $w\gamma$ is a geodesic in $G$. It is
not hard to see that the cone type of $g$ does not depend on the
choice of the geodesic $w$.

The \emph{short-lex ordering} $<$ on the set of all words $S^*$
over a finite set $S$ is a total ordering defined as follows.  For
$v,w \in S^*$, $v < w$ if $|v| < |w|$ or if $|v| = |w|$ but $v$
comes before $w$ lexicographically using some given lexicographic
order on $S$. Let $G$ be a group with finite generating set $S$.
For every $g \in G$ the smallest word $w \in S^*$ representing $g$
(with respect to the short lex ordering on $S^*$) is called the
\emph{short-lex normal form} of $g$. Note that if the set of
short-lex normal forms of $G$ is a regular language then
$\mathcal{S}_{(G,S)}$ is rational.

Let $G$ be a group with a given generating set $S$. Let $L_G
\subset S^*$ denote the set of geodesic words in $G$. If $L_G$ is
a regular language then $\mathcal{G}_{(G,S)}$ is rational. Then we
may construct a \emph{finite state automaton} (FSA) accepting
$L_G$, and from the FSA we may compute a closed form of
$\mathcal{G}_{(G,S)}$:

\begin{prop}[\cite{ECHLPT}]
\label{prop:rational_computation} Let $G$ be a group with respect
to a finite generating set $S=S^{-1}$.  Assume $L_G$ is a regular
language. Let $H$ be an FSA which accepts $L_G$. Let $M$ be the
$n\times n$ adjacency matrix of $H$ minus fail states, with the
first row of $M$ representing the adjacencies of the start state.
Then:

$$\mathcal{G}_{(G,S)}(x) = \sum_{k=0}^\infty a_kx^k = v_1(I-Mx)^{-1}v_2,$$
where $I$ is the $n\times n$ identity matrix, $v_1$ is the
$1\times n$ row vector with a one in the first column and zeros
elsewhere, and $v_2$ is the $n\times 1$ column vector of all ones.
\end{prop}

We use this theorem to compute $\mathcal{G}_{(\B,S)}$ with $S =
\{a,b\}^{\pm 1}$. Thus we first characterize geodesics for $\B$,
and then compute the appropriate FSA.

\section{The Cayley Graph of $\B$}

\subsection{Background}

Recall $\B := \langle a, b | aba = bab\rangle$. To construct the
Cayley graph for $\B$, we will need the following definitions and
theorem.

Note that in $\B$, we have that $(bab)b \equiv (aba)b = a(bab)
\equiv a(aba)$, and that $(aba)a \equiv (bab)a = b(aba) \equiv
b(bab)$.  In words, we may 'move' any occurrence of $(aba) \equiv
(bab)$ around in a word representing an element of $\B$ by
interchanging the generators $a$ and $b$. We give this property of
$\B$ a name:

\begin{defn}[The Garside property] The Garside property of $\B$ refers
to the relations $(aba)x^n \equiv  y^n(aba)$ and $(ABA)x^n \equiv
y^n(ABA)$.
\end{defn}

\begin{defn} [Right-greedy canonical form]
\cite{ECHLPT} Let $S = \{a,b,\}^{\pm 1}$. A word in $w \in S^*$ is
in \emph{right-greedy canonical form} for $\B$ if $w =
w_1w_2...w_m$ where for each $i$, $w_i \in \{a, b, ab, ba, aba,
ABA\}$, and for $i < m$:
\begin{itemize}
\item If $w_i \in \{a, ba\}$, then $w_{i+1}\in \{a, ab, aba, ABA\}$.
\item If $w_i \in \{b, ab\}$, then $w_{i+1}\in \{b, ba, aba, ABA\}$.
\item If $w_i = aba$, then $w_{i+1} = aba$.
\item If $w_i = ABA$, then $w_{i+1} = ABA$.
\end{itemize}
Denote by $RG$ the set of words in right-greedy canonical form.
\end{defn}

Note that this definition is specifically for $\B$; in
\cite{ECHLPT}, a more general definition is given.  A word $w \in
RG$ is thus of the form
$$
w = a^{k_1}b^{k_2}...a^{k_{n-1}}b^{k_n}(aba)^j,
$$
with $j \in \Z$, $k_i > 1$ for $1 < i < n$, and $k_1, k_n \geq 0$.
Every $g \in \B$ can be represented uniquely by a word in
right-greedy canonical form \cite{ECHLPT}.

\begin{prop} [\cite{DromsLewinServatius}, \cite{CrispParis}]
\label{prop:tits} The subgroup $\langle a^2,b^2 \rangle$ of $\B$
is isomorphic to the free group on two elements.
\end{prop}

This is Tits' conjecture for $B_n$ where $n = 3$. For $B_n$ with
$n \leq 5$, this result was proven in \cite{DromsLewinServatius}.
Tits' conjecture was proven in full generality in
\cite{CrispParis}. Note Proposition \ref{prop:tits} also follows
from the canonicity of the right-greedy form.

\subsection{The Graph $\Delta$}

\begin{figure}[hpbt]
\includegraphics{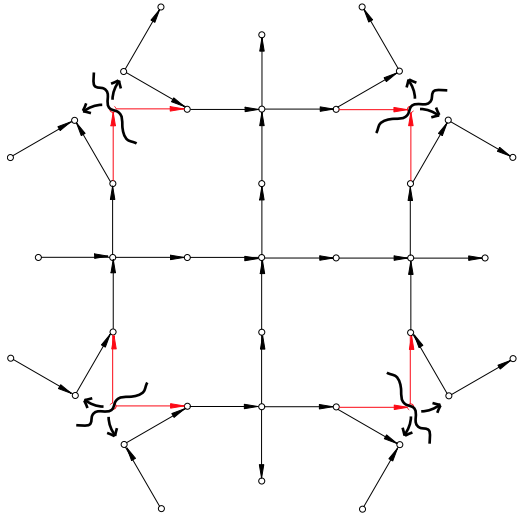}
\caption{Obtaining $C_{\B}$ from $C_{\G}$.\label{fig:C_B_3_split}}
\end{figure}

Consider the quotient group $\G$ of $\B$ obtained by adding the
relation $a^2b^2 = b^2a^2$:

$$\G = \langle a,b| aba = bab, a^2b^2 = b^2a^2\rangle.$$
Various results about this group $\G$ are proven in
\cite{Sabalka}, including a construction of the Cayley graph
$\Gamma = C_{\G}$ of $\G$, where $\Gamma$ was built by 'sewing'
together translates of the Cayley graph of the free abelian group
on $\langle a^2, b^2 \rangle$ in a certain way. The following
construction of the Cayley graph $C_{\B}$ of $\B$ is motivated by
the graph $\Gamma$, where the Cayley graph of the free abelian
group on $\langle a^2, b^2 \rangle$ is replaced by the Cayley
graph of the free group on $\langle a^2, b^2 \rangle$ (see Figure
\ref{fig:C_B_3_split}). The reader may find it enlightening to
keep $\Gamma$ in mind when reading the following description, as
$\Gamma$ may be embedded in $\R^3$, and is thus easier to picture.

\begin{defn}[The graph $T'$]
Let $T$ be the Cayley graph $C_{F_2}$ of the free group on two
generators $a$ and $b$. Let the graph $T'$ be defined to be the
subtree of $T$ spanned by $\langle a^2,b^2 \rangle$.
\end{defn}

To form the Cayley graph for $\B$, we want to take a countable
number of translates of $T'$ (one for each $k \in \Z$) and 'sew'
them together.  To do this, we define a graph $\Delta$ to have
vertex set $\bigsqcup_{k \in \Z} (VT') \times \{k\}$.  For every
$k \in \Z$, let $(ET') \times \{k\}$ denote the set of edges
induced by the inclusion $VT' \rightarrow (VT') \times \{k\}$.

A vertex $v \in T'$ is called an \emph{$a$-vertex} (respectively,
\emph{$b$-vertex}) if $v$ has exactly two incident positive edges
in $T'$, both labelled by $a$ (respectively, $b$). Define the map
$\lambda: S^* \to S^*$ to be the map induced by $\lambda(a) = b$,
$\lambda(b) = a$, $\lambda(A) = B$, $\lambda(B) = A$. Then
$\lambda^2(w) = w$. Consider any $v \in VT'$ with $w$ a path from
$1$ to $v$ in $T'$ for some $w \in S^*$. If $v$ is an $a$-vertex
(respectively, a $b$-vertex), then let $\tilde{v}$ be the vertex
represented by $\lambda(w)BA$ (respectively, $\lambda(w)AB$) in
$T'$.

\begin{figure}[hpbt]
\includegraphics{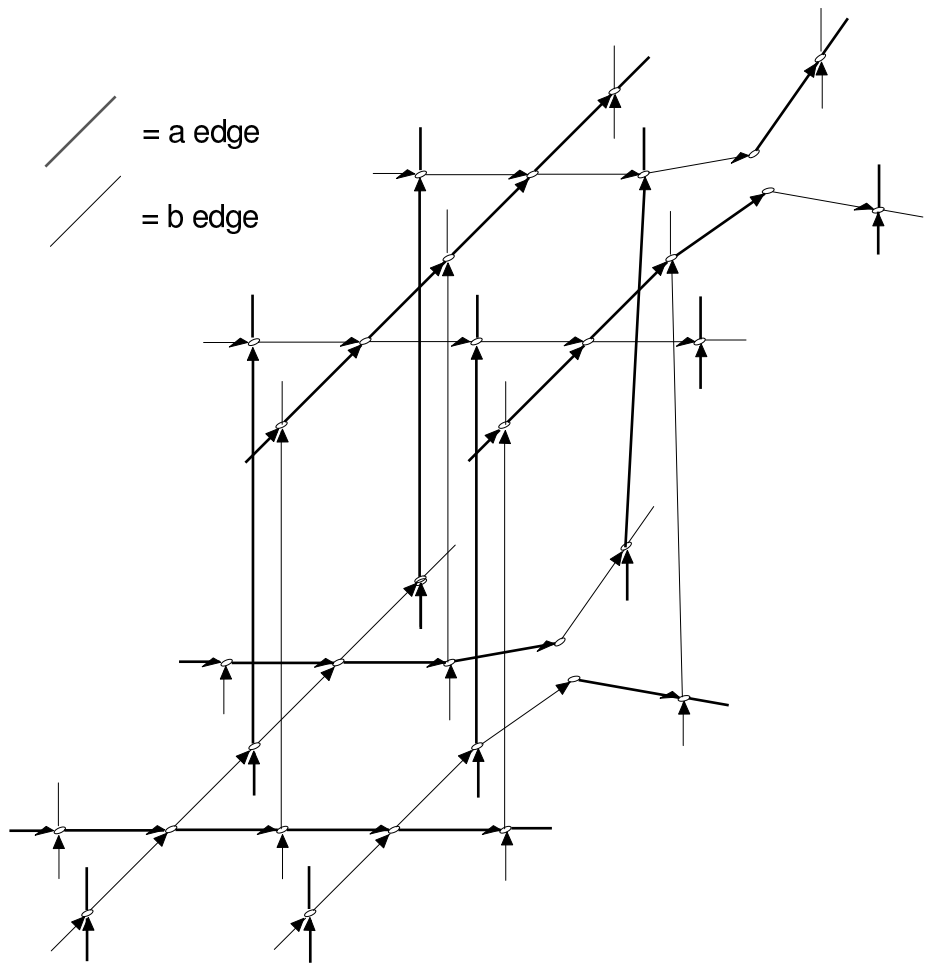}
\caption{The Cayley graph of $\B$.\label{fig:C_B_3}}
\end{figure}

\begin{defn}[The graph $\Delta$] Define the graph $\Delta$ as follows:
\begin{itemize}
\item $V\Delta := \bigsqcup_{k \in \Z} (VT') \times \{k\}$
\item The set of positively labelled edges of $\Delta$ is
 $\Bigl( \bigcup_{k \in \Z} (ET') \times \{k\} \Bigr) \cup$\\
 $\Bigl( \bigcup_{(v,k) \in V\Delta} [(v,k),(\tilde{v},k+1),a], | v$ is a~~$b$-vertex $\Bigr) \cup$\\
 $\Bigl( \bigcup_{(v,k) \in V\Delta} [(v,k),(\tilde{v},k+1),b], | v$ is an~$a$-vertex $\Bigr)$
\end{itemize}
(see Figure \ref{fig:C_B_3}).
\end{defn}

It is perhaps worth noting that this construction is related to
the fact that the pure braid group $P_3$ is isomorphic to $F_2
\times \Z$, where the $F_2$ subgroup of $P_3$ is exactly $\langle
a^2, b^2 \rangle$.

Let a word $w \in S^*$ be called \emph{almost even} if every
syllable of $w$ is even except possibly the last syllable.

\begin{prop}
\label{prop:CF} For every vertex $(v_0,k) \in V\Delta$, we may
associate to $(v_0,k)$ a unique freely reduced word $w \in S^*$
from $1$ to $(u,0)$ and whence to $(v_0,k)$ in $\Delta$ of the
form $w = w'(aba)^k$, where $w'$ is almost even and $k \in \Z$.
\end{prop}

\begin{proof}
Consider a vertex $(v_0,k) \in V\Delta$.  Let $w' \in S^*$ be a
freely reduced word such that $\overline{w'} = v_0$ in $T'$ ($w'$
is unique since $T'$ is a tree). Note that in $\Delta$, the path
labelled $aba$ beginning at $(v_0,k)$ ends at $(v_0',k+1)$, where
$v_0'$ is the vertex represented by $\lambda(w')$ in $T'$. Let $u$
be the vertex in $VT'$ represented by $\lambda^k(w')$ (i.e. $u =
v_0$ if $k$ is even, and $u = v_0'$ if $k$ is odd). Consider the
path $w := \lambda^k(w')(aba)^k$ in $\Delta$. The path $w$ is a
path from $(1,0)$ to $(u,0)$ and whence to $(v_0,k)$, where the
first part of $w$ - $\lambda^k(w')$ - is also a path in $T'$.  In
$\Delta$, this means $\lambda^k(w')$ only connects vertices in
$\Delta$ of the form $(v,0)$. To stay in $T'$, $\lambda^k(w')$
must be almost even, by the definition of $T'$.
\end{proof}

Let $CF$ (for `Cartesian Form') denote the set of words $w$ with
$w = w'(aba)^k$ where $w'$ has only even syllables except possibly
the last syllable, and $k \in \Z$.

\begin{thm}
\label{thm:Delta_is_C_B_3} The graph $\Delta$ is the Cayley graph
for $\B$ with respect to the generating set $S = \{a,b\}^{\pm 1}$.
\end{thm}

\begin{proof}
By Proposition \ref{prop:CF}, vertices of $\Delta$ are in
bijective correspondence with elements of $CF$, and elements of
$RG$ are in bijective correspondence with vertices of the Cayley
graph of $\B$. Thus, to construct a bijection from $V\Delta$ to
$VC_{\B}$, we will define a bijection between $CF$ and $RG$.

We define a map $\Phi_1: CF \mapsto RG$ as the result of the
following rewriting algorithm.  Consider $w \in CF$, where $w =
w'(aba)^k$, with $w'$ almost even and $k \in \Z$.

First, this algorithm deals with every occurrence of $AB$ or $BA$
in $w'$. These occurrences are non-overlapping since $w'$ is
almost even. For each occurrence of $AB$ in $w'$, insert $bB$
beforehand, creating a power of $ABA$: $AB \cong (bB)(AB) = b(BAB)
\equiv b(ABA)$. Similarly, for every occurrence of $BA$ insert
$aA$ beforehand, again creating a power of $ABA$. Then, pull all
created occurrences of $ABA$ to the right end of $w$ using the
Garside property. The result of the algorithm thus far is a word
of the form $w''(aba)^{k'}$, where $w''$ contains no occurrences
of $AB$ or $BA$; in other words, all negative syllables are
isolated in $w''$. Note if $s^c$ ($s \in \{a,b\},c > 0$) is a
syllable in $w''$ with $c = 1$, then $s^c$ is the result of one of
the insertions just described, and thus $s^c$ is between two
negative syllables. For instance, $AABBBBAA \mapsto
A(bB)ABBB(aA)BAA \mapsto AbA^2bA(ABA)^2$.

The second part of the algorithm deals with any remaining negative
syllables in $w''$ as follows. For every occurrence of $A$, insert
the trivial words $bB$ and $Bb$ before and after the $A$,
respectively, thus creating an occurrence of $BAB \equiv ABA$.
Similarly, for every occurrence of $B$ insert $aA$ beforehand and
$Aa$ after. Finally, pull all powers of $ABA$ to the end of $w''$
using the Garside property. For instance, consider $w'' = A^c$ for
$c$ odd.  We have $w'' = A^c \equiv (bBABb)(bBABb)\dots (bBABb)$,
where there are $c$ occurrences of $(bBABb)$.  Pulling each power
of $BAB \equiv ABA$ to the end, we get $ba^2b^2\dots b^2a(ABA)^k$.
Note that for each negative syllable, this process creates
occurrences of $ABA$, $a^2$, and $b^2$, and adds 1 to the power of
adjacent syllables. Increasing the power of adjacent syllables
eliminates the occurrences of syllables $s^c$ with $c = 1$ in
$w''$ created by the first part of the algorithm.

So, $w$ is equal in $\B$ to the resultant word of the form
$\tilde{w} := w'''(aba)^{c''}$, where $w'''$ has only positive
syllables of exponent strictly greater than 1 (except perhaps the
first or last syllables). Thus, $\tilde{w} \in RG$, and we set
$\Phi_1(w) = \tilde{w}$.

We now define the map $\Phi_2: RG \mapsto CF$ by the following
rewriting algorithm. Consider $\tilde{w} \in RG$, where $\tilde{w}
= w'(aba)^c = a^{k_1}b^{k_2}...a^{k_{n-1}}b^{k_n}(aba)^c$ with $j
\in \Z$, $k_i > 1$ for $1 < i < n$, and $k_1, k_n \geq 0$. If
$k_i$ is even for each $1 \leq i < n$, then $\tilde{w}$ is already
in $CF$, and we set $\Phi_2(\tilde{w}) = \tilde{w}$. Otherwise, we
proceed as follows.

Let $j = j_w$ be the least $i$ such that $k_i$ is odd in $w'$.
Without loss of generality assume the $j$th syllable is $a^{k_j}$.
If $k_{j+1}$ is even, then:

\begin{eqnarray*}
w&\equiv &\dots a^{k_j}b^{k_{j+1}} \dots \\
 &\equiv &\dots a^{k_j-1}(BBbb)abb^{k_{j+1}-1} \dots \\
 &\equiv &\dots a^{k_j-1}B^2ba^{k_{j+1}-1}(bab) \dots \\
 &\equiv &\dots a^{k_j-1}B^2A^2aa^{k_{j+1}-2}(bab)^2 \dots \\
 &\equiv &\dots a^{k_j-1}B^2A^2B^2bb^{k_{j+1}-3}(bab)^3 \dots \\
 &\equiv &\dots\\
 &\equiv &\dots a^{k_j-1}B^2A^2B^2\dots A^2B(aba)^{k_{j+1}}\dots \\
\end{eqnarray*}
where there are $k_{j+1}-2$ negative squared syllables. If
$k_{j+1}$ is odd, then the result of the same process is $\dots
a^{k_j-1}B^2A^2B^2\dots B^2A(aba)^{k_{j+1}}\dots$. Let $w''$
denote this result with the powers of $aba$ pulled to the right
end, using the Garside property.  Then $\tilde{w} \equiv w''$. If
$w''$ has $n'$ syllables followed by $(aba)^{c+k_{j+1}}$, then
note $n - j_{w'} > n' - j_{w''}$. Thus this process will terminate
upon iteration.

Continue to apply this algorithm to each resultant word until no
odd syllables are left (except possibly the last one). The result
of this process is a word of the form $w(aba)^c$ where $w$ is
almost even and $k \in \Z$, as desired. Set $\Phi_2(\tilde{w}) =
w(aba)^k$.

\begin{lemma}
\label{lem:bijective_correspondence} The functions $\Phi_1: CF \to
RG$ and $\Phi_2: RG \to CF$ are bijections, and $\Phi_2 =
\Phi_1^{-1}$.
\end{lemma}

\begin{proof}
We have already defined maps $\Phi_1: CF \mapsto RG$ and $\Phi_2:
RG \mapsto CF$. We need to show that $\Phi_1 \circ
\Phi_2(\tilde{w}) = \tilde{w}$ for every $\tilde{w} \in RG$ and
$\Phi_2 \circ \Phi_1(w) = w$ for every $w \in CF$. Both $\Phi_1$
and $\Phi_2$ only use the relations of $\B$; therefore, for any
$\tilde{w} \in RG$, $\Phi_1 \circ \Phi_2 (\tilde{w}) \equiv
\tilde{w}$. Since $RG$ is a set of canonical normal forms for
$\B$, this implies $\Phi_1 \circ \Phi_2 (\tilde{w}) = \tilde{w}$.
Thus it suffices to show $\Phi_2 \circ \Phi_1(w) = w$ for every $w
\in CF$.

Consider $\Phi_2 \circ \Phi_1(w)$ for some non-trivial $w =
w'(aba)^k \in CF$. Then $w$ is of the form
$w^+_1w^-_1w^+_2w^-_2\dots w^+_mw^-_m$ for some $m \geq 0$, where
each $w^+_k$ (respectively, $w^-_k$) is a word with only positive
(respectively, negative) syllables.  We prove $\Phi_2 \circ
\Phi_1(w) = w$ by induction on $m$. Note if $m = 0$ or $m = 1$ and
$w^-_1 = 1$ there is nothing to prove.

Assume that for $m = 1,\dots, n-1$ we have $\Phi_2 \circ \Phi_1(w)
= w$.  For $m = n$, we have $w = w^+_1w^-_1w_1$ where $w_1 \in CF$
and $\Phi_2 \circ \Phi_1(w_1) = w_1$.  Note $$\Phi_2 \circ \Phi_1
(w) = \Phi_2(w^+_1\Phi_1(w^-_1w_1)) = w^+_1[\Phi_2 \circ
\Phi_1](w^-_1w_1),$$ and $\Phi_1(w^-_1w_1) =
\Phi_1(w^-_1)\Phi_1(w_2)$, where $w_2 = w_1$ or $w_2 \equiv
\lambda(w_1)$. Thus, since $\Phi_2$ scans from left to right, if
$\Phi_2 \circ \Phi_1 (w^-_1) = w^-_1$, then $\Phi_2 \circ \Phi_1
(w_2) = \Phi_2 \circ \Phi_1 (w_1) = w_1$, and $\Phi_2 \circ
\Phi_1(w) = w$.  So it suffices to show that $\Phi_2 \circ \Phi_1
(w^-_1) = w^-_1$.

We will prove $\Phi_2 \circ \Phi_1 (w^-_1) = w^-_1$ by induction
on the number of negative syllables in $w^-_1$. Consider the
number of (negative) syllables of $w^-_1$. Assume without loss of
generality that the first syllable of $w^-_1$ is $A^{2k_1}$, $k_1
> 0$. If $w^-_1$ has only one syllable, then $\Phi_1(w^-_1) =
ba^2b^2a^2...a^2b(ABA)^{2k_1}$, where there are $k_1$ occurrences
of $a^2$, and:

\begin{eqnarray*}
\Phi_2(\Phi_1(w^-_1))&=&\Phi_2(ba^2b^2a^2\dots a^2b(ABA)^{2k_1})\\
                     &=&A^2\Phi_2(aba^2b^2a^2\dots b^2a(ABA)^{2k_1-1})\\
                     &=&A^2\Phi_2(B^2bb^2a^2b^2\dots a^2b(ABA)^{2k_1-2})\\
                     &=&A^2\Phi_2(ba^2b^2\dots a^2b(ABA)^{2k_1-2})\\
                     &=&\dots\\
                     &=&A^{2k_1}\\
\end{eqnarray*}
as required.

If $w^-_1$ has two negative syllables, it is of the form
$A^{2k_1}B^{2k_2}$. Then the first part of the algorithm defining
$\Phi_1$ eliminates the occurrence of $AB$, yielding
$A^{2k_1}bA^{2k_2}(ABA)$. Then, $\Phi_1(w^-_1) = ba^2b^2\dots
b^2a^3b^2\dots a^2b(ABA)^{2k_1+2k_2-1}$, where there are $k_1-1$
occurrences of $a^2b^2$ before the $a^3$ and $k_2-1$ occurrences
of $b^2a^2$ after. Thus:

\begin{eqnarray*}
\Phi_2(\Phi_1(w^-_1))&=&\Phi_2(ba^2b^2a^2\dots b^2a^3b^2a^2\dots a^2b(ABA)^{2k_1+2k_2-1})\\
                     &=&A^2\Phi_2(ba^2\dots b^2a^3b^2a^2\dots a^2b(ABA)^{2k_1+2k_2-1-1})\\
                     &=&\dots\\
                     &=&A^{2k_1-2}\Phi_2(ba^3b^2\dots a^2b(ABA)^{2k_2-1-(2k_1-2)})\\\\
                     &=&A^{2k_1-2}A^2\Phi_2(a(aba)aab^2\dots a^2b(ABA)^{2k_2+1})\\\\
                     &=&A^{2k_1}\Phi_2(ab^2a^2\dots b^2a(ABA)^{2k_2})\\\\
                     &=&A^{2k_1}B^{2k_2},\\
\end{eqnarray*}
again as required.

If there are more than two negative syllables in $w^-_1$, then
$w^-_1$ is of the form $A^{2k_1}B^{2k_2}A^{2k_3}...$. Then,
$\Phi_1(w^-_1) = ba^2b^2\dots b^2a^3b^2a^2 \dots
b^2a^3b^2a^2\dots$, with $k_j-1$ occurrences of $b^2$ between the
$(j-1)$th and the $j$th occurrence of $a^3$.  When $\Phi_2$ is
applied to this result, scanning from left to right we get

\begin{eqnarray*}
\Phi_2(\Phi_1(w^-_1))&=&\Phi_2(\Phi_1(A^{2k_1}B^{2k_2}A^{2k_3}...))\\
                     &=&\Phi_2(ba^2b^2\dots b^2a^3b^2a^2 \dots b^2a^3b^2a^2\dots)\\
                     &=&A^{2k_1}\Phi_2(ab^2a^2 \dots b^2a^3b^2a^2\dots)\\
                     &=&A^{2k_1}\Phi_2(\Phi_1(B^{2k_2}A^{2k_3}...)),\\
\end{eqnarray*}
and we may induct on the number of negative syllables of $w^-_1$.
Thus $\Phi_2(\Phi_1(w^-_1))$ preserves $w^-_1$, and Lemma
\ref{lem:bijective_correspondence} is proved.
\end{proof}

We now return to the proof of Theorem \ref{thm:Delta_is_C_B_3}. It
is left to prove there is a bijection between edges which
preserves adjacency, i.e. that for $w\in CF$, if $\Phi_1(w) = u
\in RG$, then $wa$ is mapped to $ua$, and $wb$ is mapped to $ub$.
We claim that indeed $\Phi_1$ and $\Phi_2$ induce a graph
isomorphism between $\Delta$ and $C_{\B}$.

Let $w\in CF$ be mapped to $\Phi_1(w) = u \in RG$. For a given $x
\in \{a,b\}^{\pm 1}$, let $w_1 \in CF$ and $u_1 \in RG$ be such
that $w_1 \equiv wx$ and $u_1 \equiv ux$. We want $\Phi_1(w_1) =
u_1$ if and only if $\Phi_1(w)x \equiv \Phi_1(w_1)$.

There are many cases to deal with for this argument. We deal with
two; the rest may be argued similarly. Assume without loss of
generality that $x = a$ and $w$ is of the form $w = \dots
b^{k_{n-1}}a^{k_n}(aba)^{i}$, with $k_i$ even for $1 \leq i < n$.
Look at $\Phi_1(w_1)$. Assume $k_n < 0$, $k_n$ even, $k_{n-1} >
0$, and $j_u$ even. Since $k_n < 0$, $u = \dots a^2b(aba)^{j}$.
Then:

\begin{itemize}
\item If $j_w$ is even, then, in $\B$, $wa \equiv \dots
a^{k_n}(aba)^{j_w}a \equiv \dots a^{k_n+1}(aba)^{j_w} \in CF$, so
$w_1$ is of the form $\dots a^{k_n+1}(aba)^{j_w}$. But

$$\Phi_1(w_1) \equiv \dots a(aba)^{j_u+1} \equiv \dots
a^2ba(aba)^{j_u} \equiv a^2b(aba)^{j_u}a \equiv ua = \Phi_1(w)a$$
by definition, and so $u_1 = \Phi_1(w_1)$ as required.

\item If $j_w$ is odd, then, in $\B$, $wa \equiv \dots
a^{k_n}(aba)^{j_w}a \equiv \dots a^{k_n}b(aba)^{j_w} \in CF$ since
$k_n$ is even. Thus $w_1$ is of the form $\dots
a^{k_n}b(aba)^{j_w}$.  Again,

$$\Phi_1(w_1) \equiv \dots a^2b^2(aba)^{j_u} \equiv \dots
a^2b(aba)^{j_u}a \equiv ua = \Phi_1(w)a,$$ and so $u_1 =
\Phi_1(w_1)$ as required.
\end{itemize}

In every case we see $\Phi_1$ preserves vertex adjacency, making a
bijection between the two graphs which preserves adjacency. Thus,
$\Delta \cong C_{\B}$, and the proof of Theorem
\ref{thm:Delta_is_C_B_3} is complete.
\end{proof}

\section{Characterizing Geodesics of $\B$}

\subsection{Some Definitions}

Looking at the Cayley graph $\Delta$ of $\B$, one might think that
a word $w = w'(aba)^k \in CF$ (where $w'$ is almost even and $k
\in \Z$) is geodesic. However, this is not necessarily the case.
For instance, the word $a^2b^2A^2B^2 \in CF$ is equal in $\B$ to
the shorter word $aBaBaB$.

It is sometimes shorter in $\Delta$ to move from the translate
$T'\times \{k\}$ of $T'$ to a nearby translate $T'\times \{k \pm
1\}$ and then back. Algebraically, we see that if $w$ is of the
form $w = ...(ab)b^{k_1}a^{k_2}...a^{k_n}(AB)...$, then,

\begin{eqnarray*}
w&=&...(ab)b^{k_1}a^{k_2}...a^{k_n}(AB)...\\
 &\equiv &...B(bab)b^{k_1}a^{k_2}...a^{k_n}(ABA)a...\\
 &\equiv & ...Ba^{k_1}b^{k_2}...b^{k_n}a...\\
\end{eqnarray*}
by the Garside property.  This gives a shorter representation, so
if $w$ contains both $ab$ and $AB$, then $w$ is not geodesic. In
fact, if $w$ contains both one of $\{ab,ba\}$ and one of
$\{AB,BA\}$ as subwords, then a similar argument proves $w$ is not
geodesic.

\begin{defn}[$*^+$ and $*^-$]
Let $*^+$ $:= \{ab, ba\}$ and $*^-$ $:= \{AB, BA\}$.
\end{defn}

Thus a geodesic may not have as subwords elements of both $*^+$
and of $*^-$. Indeed, also by the Garside property, a geodesic may
not have as subwords both $(aba)$ and either $A$ or $B$, and it
may not have as subwords both $(ABA)$ and either $a$ or $b$. This
motivates the following:

\begin{defn}[The * condition]
\label{def:*} An element $w \in S^*$ is said to satisfy the
\emph{* condition} if $w$ does not contain as subwords elements of
both $*^+$ and $*^-$.
\end{defn}

\begin{defn}[The ** condition]
\label{def:**} An element $w \in S^*$ is said to satisfy the
\emph{** condition} if:
\begin{itemize}
\item $w$ does not contain as subwords both $(aba)$ and either $A$ or $B$, and
\item $w$ does not contain as subwords both $(ABA)$ and either $a$ or $b$.
\end{itemize}
\end{defn}

\begin{defn}[$\Geods$]
We define $\Geods$ to be the set of all $w \in S^*$ such that $w$
is freely reduced and satisfies both the
* condition and the ** condition.
\end{defn}

The following is now obvious:

\begin{thm}
We have that $L_{\B} \subseteq \Geods$.
\label{thm:geodesicspartial}
\end{thm}

The * and ** conditions thus severely limit the set of possible
geodesics. They are necessary conditions on a geodesic; we will
eventually see that they are also sufficient.

Recall that in Section 2 we discussed the short-lex ordering on
$S^*$. Fix the lexicographical ordering $a<A<b<B$ on $S =
\{a,b\}^{\pm 1}$.  Let $SL$ denote the set of short-lex normal
forms for $\B$.

\begin{lemma}
\label{lem:SL} A word $w$ is in $SL$ if and only if it is in one
of the following forms:
\begin{itemize}
\item $(a^i)(b^{j_1}A^{j_2}b^{j_3}...A^{j_{m-1}}b^{j_m})(a^{k_1}b^{k_2}a^{k_3}...a^{k_n}b^k)$,
\item $(a^i)(b^{j_1}A^{j_2}b^{j_3}...A^{j_{m-1}}b^{j_m})(a^{k_1}b^{k_2}a^{k_3}...b^{k_n}a^k)$,
\item $(a^i)(B^{j_1}a^{j_2}B^{j_3}...a^{j_{m-1}}B^{j_m})(A^{k_1}B^{k_2}A^{k_3}...A^{k_n}B^k)$, or
\item $(a^i)(B^{j_1}a^{j_2}B^{j_3}...a^{j_{m-1}}B^{j_m})(A^{k_1}B^{k_2}A^{k_3}...B^{k_n}A^k)$,
\end{itemize}
where $i \in \Z$, $m$ and $n$ are non-negative integers, $j_l > 0$
for every $0 \leq l \leq m$, $k_l > 1$ for every $0 \leq l \leq
n$, and $k$ is either $0$ or $1$.
\end{lemma}

\begin{proof}
These are computed from the FSA for $SL$, which was computed by
the program KBMAG.
\end{proof}

Note that $SL$ is a subset of the set of geodesics of $\B$.

\subsection{Proving $\Geods = L_{\B}$}

Recall that the set of geodesic words of $\B$ is $L_{\B}$.  With
the definitions above, we may now restate Theorem
\ref{thm:geodesics} with simpler notation:

\begin{thm}
We have that $\Geods = L_{\B}$. In other words, a freely reduced
word $w \in S^*$ is geodesic in $\B$ if and only if $w$ satisfies
the * and ** conditions.
\end{thm}

\begin{proof}
By Theorem \ref{thm:geodesicspartial}, we have that $L_{\B}
\subseteq \Geods$.  It remains to prove that $\Geods \subseteq
L_{\B}$.  We do so by showing that a word $w$ in $\Geods$ is equal
in $\B$ to a word $w'$ in $SL$ obtained from $w$ by a
length-preserving map. We present an algorithm for rewriting $w$
to $w'$ in two steps. The first step takes a word $w \in \Geods$
into a temporary form (increasing the length by a certain amount),
and the second step maps words in the temporary form to $SL$
(decreasing the length by the same amount).

Let $TF$ denote the set of all words of the form
$(x^{k_1}Y^{k_2}...Y^{k_{n-1}}x^{k_n})(aba)^j$, with $k_n \geq 0$,
$k_i \geq 1$ for $1 \leq i < n$, and for some
$(x,y)\in\{(a,b),(b,a)\}^{\pm 1}$.  We will say that a word in
$TF$ is in \emph{temporary form}.

Define $\Psi_1: S^* \rightarrow TF$ by the following algorithm.
For a word $w \in S^*$, first freely reduce $w$.  Then, pull all
subwords of the form $aba$, $bab$, $ABA$, and $BAB$ to the right
end using the Garside property, so $w \equiv w'(aba)^{j'}$, and
freely reduce $w'$. Then, look at $w'$ in terms of its syllables.
Working from left to right, if there are ever two syllables in a
row of the form $x^{k_i}y^{k_{i+1}}$ with $(x,y) \in
\{(a,b),(b,a)\}^{\pm 1}$ , then since

\begin{eqnarray*}
x^{k_i}y^{k_{i+1}}&=&x^{k_{i}-1}xyy^{k_{i+1}-1}\\
&\equiv &x^{k_{i}-1}Y(yxy)y^{k_{i+1}-1}\\
&\equiv &x^{k_{i}-1}Yx^{k_{i+1}-1}(yxy),
\end{eqnarray*}
we may replace $x^{k_i}y^{k_{i+1}}$ with
$x^{k_i-1}Yx^{k_{i+1}-1}(yxy)$ and pull the power of $(aba)$ out
of $w'$ to the end of $w$. Then repeat this replacement for every
appropriate pair of syllables, working to the end of $w'$. The
resulting word is $\Psi_1(w)$.

There are a few things to note about $\Psi_1$ concerning the
length in $\B$ of $\Psi_1(w)$.  The initial process of moving
occurrences of $(xyx)$ to the end of the word does not change the
length of a freely reduced word $w$, unless there are both
positive and negative powers of $aba$ in the original word.  For
$w \in \Geods$, however, this does not happen by the ** condition.
For the remainder of the algorithm, each described replacement
increases by two the length of the resultant word, \emph{unless}
cancellation occurs. Cancellation may occur in two locations.  If
the power of $aba$ pulled to the end has opposite sign of another
power of $aba$ also pulled to the end, then cancellation occurs.
However, this only happens when initially $w$ had as subwords
elements of both $*^+$ and $*^-$. This contradicts the
* condition, so for $w \in \Geods$, cancellation does not occur at
the end of the word. Cancellation could also occur if $k_i$ or
$k_{i+1}$ equals $1$. If $k_i = 1$ and the $i$th syllable is $x$,
and if the $(i-1)$th syllable of $w$ exists and is $y^{k_{i-1}}$
with $k_{i-1} > 0$, then $Y$ will cancel with the $(i-1)$th
syllable. Similarly, if $k_{i+1} = 1$ and if the $(i+2)$th
syllable exists and is a positive power of $x$, cancellation will
again occur. However, in each case $w'$ would have to have a
subword of the form $yxy$ or $xyx$, respectively, and by the first
part of the algorithm, this does not happen. Thus, for $w \in
\Geods$, the length of $\Psi_1(w)$ is exactly the length of $w$
plus twice the number of occurrences of subwords of the form $ab$,
$ba$, $AB$, and $BA$. Let $\lambda$ denote the number of such
subwords.

We define $\Psi_2: TF \rightarrow SL$ as the result of the
following rewriting algorithm.  Consider a word $w' \in TF$.  Then
$w'$ is of the form

$$(x^{k_1}Y^{k_2}...Y^{k_{n-1}}x^{k_n})(aba)^{j}$$
for some $(x,y) \in \{(a,b),(b,a)\}^{\pm 1}$. Note that
$x^{k_1}Y^{k_2}...Y^{k_{n-1}}x^{k_n} \in SL$, so if $j = 0$,
define $\Psi_2(w') = w'$.  For $j \neq 0$, without loss of
generality assume $j$ is positive and assume $x = a$ (for the
other seven combinations, similar definitions of $\Psi_2$ may be
made, but are not elaborated upon here). Then $w'$ is of the form
$(a^{k_1}B^{k_2}...B^{k_{n-1}}a^{k_n})(aba)^{j}$. It is still true
that $a^{k_1}B^{k_2}...B^{k_{n-1}}a^{k_n} \in SL$; thus, we need
to deal with $(aba)^{j}$.  First, move one of the powers of
$(aba)$ \emph{before} the \emph{first} negative syllable and
freely reduce (if no such syllable exists, skip this step):

\begin{eqnarray*}
w'&=&(a^{k_1}B^{k_2}a^{k_3}...a^{k_{n-2}}B^{k_{n-1}}a^{k_n})(aba)^{j}\\
&\equiv &(a^{k_1}(aba)A^{k_2}b^{k_3}...b^{k_{n-2}}A^{k_{n-1}}b^{k_n})(aba)^{j-1}\\
&\equiv &(a^{k_1+1}bA^{k_2-1}b^{k_3}...b^{k_{n-2}}A^{k_{n-1}}b^{k_n})(aba)^{j-1}\\
&\equiv&(a^{k_1+1})(bA^{k_2-1}b^{k_3}...b^{k_{n-2}}A^{k_{n-1}}b^{k_n})(aba)^{j-1}.
\end{eqnarray*}
Again,
$(a^{k_1+1})(bA^{k_2-1}b^{k_3}...b^{k_{n-2}}A^{k_{n-1}}b^{k_n})
\in SL$, and the length of the word is decreased by two.

If $j > 1$, then the remaining powers of $(aba)$ will be
eliminated as follows.  While there are still powers of $(aba)$
left, if there are any negative syllables remaining, take an
$(aba)$ and place it \emph{after} the \emph{last} negative
syllable using the Garside property and then freely reduce.  For
instance, for $w = (a^2B^2aB^2)(aba)^2$, initially we make the
word $(a^2(aba)A^2bA^2)(aba) = (a^3)(bAbA^2)(aba)$, and then
$(a^3)(bAbA^2(aba)) = (a^3)(bAbAb)(a)$ (which is in short-lex
form). We may continue to reduce the exponent on $(aba)$ in this
process, each time decreasing the length of the word by two, until
either there are no more $(aba)$s, or there are no more negative
syllables remaining. Let $\alpha$ denote the number of occurrences
of $A$ and $B$ originally in $w$. Then thus far we have decreased
the length of the word by $min\{\alpha,j\}$.

If $j > \alpha$, to eliminate the remaining powers of $(aba)$,
insert them one by one \emph{after} the first syllable. In
short-lex normal form, only the second syllable may be raised to
the first power (see Definition \ref{lem:SL}).  If we begin with a
word in $SL$, though, this algorithm preserves that property, and
the result is in $SL$. In total, then, $\Psi_2$ shortens a word
$w' \in TF$ by $2*min\{\alpha,j\}$.

Now consider the map $\Psi := \Psi_2\circ\Psi_1$.  Let $w_1 :=
\Psi_1(w)$ and $w_2 := \Psi_2(w_1) = \Psi(w)$. When restricted to
$w \in \Geods$, as already noted for $w_1$, $\Psi_1$ first pulls
out all powers of $(aba)$ (preserving the length), and then pulls
out all subwords of the form $ab$, $ba$, $AB$, or $BA$, each time
increasing the length of the word by two. For a word in $\Geods$,
if there exist any subwords of the form $aba$ or $bab$, then by
the ** condition there are no negative powers in $w$, and $j \geq
\alpha = \lambda$, so

$$
|w_2| = |w_1| - 2*min(\alpha,j) = (|w| + 2*\lambda) - 2*\alpha =
|w|
$$
and in this case $\Psi$ preserves the length of $w$.

If $w$ does not have any subwords of the form $aba$ or $bab$, then
either $w$ is already in $SL$ (as shown above) or $w$ contains
subwords from exactly one of $*^+$ or $*^-$. Without loss of
generality, assume $w$ contains elements from only $*^+$.  Then $j
= \lambda$ (recall $j$ the exponent of $(aba)$ in $w_1 =
\Psi(w)$). Furthermore, for every occurrence of an element of
$*^+$ in $w$, $\Psi_1$ inserts exactly one negative letter which
does not cancel with its neighbors. Thus we have $\alpha \geq j$
(recall $\alpha$ is the number of negative elements in $w_1$), and
$|w_2| = |w_1| - 2*min\{\alpha,j\} = (|w|+2*\lambda)-2*j = |w|$,
as desired.

Thus, for any $w \in \Geods$, $|\Psi(w)| = |w|$, and $\Psi$
restricted to $\Geods$ is a length-preserving map, and we have
proven Theorem \ref{thm:geodesics}.
\end{proof}

\subsection{Translation lengths in $\B$}

For a general group $G$ with generating set $S$, we define the
length $|g|_S$ of an element $g \in G$ with respect to $S$:

$$|g|_S := \inf_{w\in S^*, \overline{w} = g}\{|w|\}.$$
The \emph{translation length} of an element $g\in G$ is defined to
be:

$$\tau_S(g) := \limsup_{n\to \infty} \frac{|g^n|_S}{n}.$$

It is well known that the limit here always exists, by
subadditivity. It is straightforward to check that that
$\tau_S(g)$ depends only on the conjugacy class $CC(g)$ of $g$,
where $CC(g) := \{h^{-1}gh | h \in G\}$. For more information on
translation numbers, see for instance \cite{GerstenShort},
\cite{Conner}, or \cite{Kapovich}.

It is known that for any word-hyperbolic group $G$, for any finite
generating set $S$ of $G$ and any $g \in G$, $\tau_S(g)$ is a
rational number \cite{BaumslagGerstenShapiroShort}; the
denominator of $\tau_S(g)$ is bounded (stated by Gromov; proved in
\cite{Swenson}).  The previous results of this section now also
allow us to prove:

\begin{corollary} For $S = \{a,b\}^{\pm 1}$ and any element $g \in \B$,
$\tau_S(g)$ is an integer. In particular, $\tau_S(g) = |g'|_S$,
where $g'$ is a shortest element in $CC(g)$.
\end{corollary}

\begin{proof}
Consider a general element $g \in \B$.  We may assume that $g$ is
the smallest element in its conjugacy class.  Then a word $w \in
S^*$ representing $g$ is geodesic if and only if $w$ is cyclically
reduced, any cyclic permutation of $w$ is geodesic, and $w$
satisfies the * and ** conditions.  For any $n \in \N$, if $w^n$
is not geodesic, then $w^n$ violates either the * or the **
condition.  By the periodicity of $w^n$, there exists some subword
$w'$ of $w^n$ of length $|w|$ which violates either the * or the
** condition. But then $w'$ is a cyclic permutation of $w$ which
is not geodesic, contradicting the choice of $g$. Thus $w^n$ is a
geodesic representative of $g^n$, and so $|g^n|_S = n|w|$ for each
$n$, and $\tau_S(g) = |w| = |g|_S$.
\end{proof}

\section{Computing the Growth Series}

\subsection{The FSA for Geodesics of $\B$}

\begin{figure}[hpbt]
\includegraphics{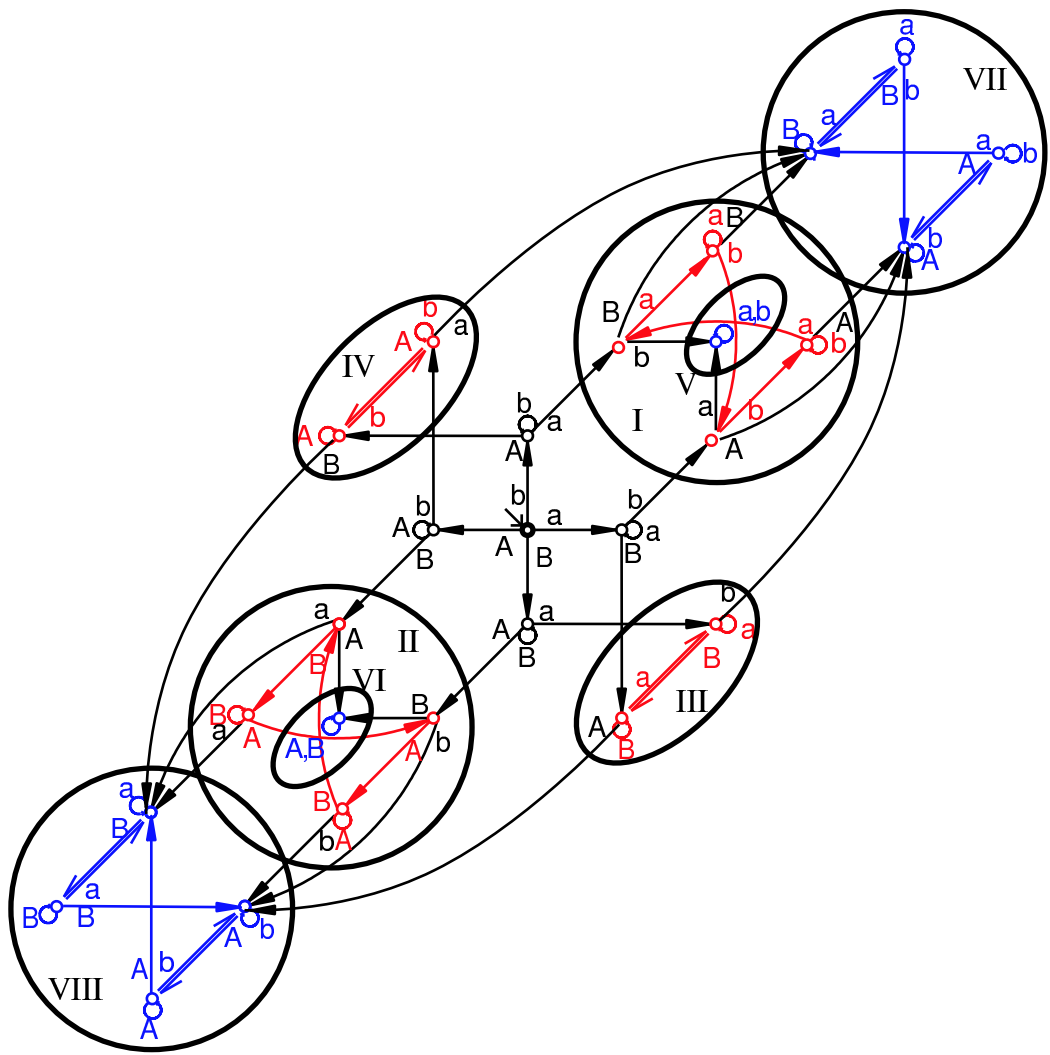}
\caption{The FSA which accepts geodesics of $\B$, with its
different components labelled.  All states shown are accept
states; the single fail state is not shown.  There are 27 states.
\label{fig:FSA_B_3}}
\end{figure}

By Theorem \ref{thm:geodesics} we know that a freely reduced word
$w \in S^*$ is a geodesic in $\B$ if and only if $w$ satisfies the
* and ** conditions. Using this knowledge, constructing the FSA
accepting geodesic words for $\B$ is straightforward - the FSA
simply needs to keep track of which choices have been made with
respect to the * and ** condition. It is clear that checking
whether a word satisfies the * and ** conditions requires keeping
track of only a finite amount of information, so such an FSA
exists. We have:

\begin{thm}
The FSA in Figure \ref{fig:FSA_B_3} accepts a word $w \in S^*$ if
and only if $w$ is a geodesic of $\B$. Thus the set of geodesics
$L_{\B}$ of $\B$ is regular and $\mathcal{G}_{(\B,S)}$ is
rational.
\end{thm}

The different parts of Figure \ref{fig:FSA_B_3} represent which
decisions have been made by the automaton when reading a word $w
\in \Geods$ with respect to the * and ** condition. The decisions
made are described below:

\begin{itemize}
\item Part I:  $w$ has as a subword an element of $*^+$.
\item Part II: $w$ has as a subword an element of $*^-$.
\item Part III: $w$ has as subwords $a$ and
$B$, but no elements of $*^+$ or $*^-$.
\item Part IV: $w$ has as subwords $b$ and
$A$, but no elements of $*^+$ or $*^-$.
\item Part V: $w$ has as a subword either $aba$ or $bab$.
\item Part VI: $w$ has as a subword either $ABA$ or $BAB$.
\item Part VII: $w$ has as subwords an element of $*^+$
and a negative syllable.
\item Part VIII: $w$ has as subwords an element of $*^-$
and a positive syllable.
\end{itemize}

\subsection{The Geodesic Growth Series of $\B$}

Recall that by Proposition \ref{prop:rational_computation}, we
have

$$\mathcal{G}_{(\G,S)}(x) = v_1(I-Mx)^{-1} v_2$$
with $I$, $M$, $v_1$ and $v_2$ as in Proposition
\ref{prop:rational_computation}.

When the adjacency matrix of the FSA for $\B$ is plugged into this
formula using Maple 4.5, we get:

\begin{thm} (c.f. Theorem \ref{thm:GGS_B})
The Geodesic Growth Series of $\B$ with respect to $S =
\{a,b\}^{\pm 1}$ is

$$\mathcal{G}_{(\B,S)}(x) = \frac{x^4+3x^3+x+1}{(x^2+x-1)(x^2+2x-1)}.$$
\end{thm}

For purposes of comparison, using the method proved in
\cite{Brazil} and using the FSA which accepts the short-lex normal
form of $\B$ computed by the program KBMAG (see
\cite{EpsteinHoltRees}), we have:

\begin{thm}
The Spherical Growth Series of $\B$ is

$$\mathcal{S}_{(\B,S)}(x) = \frac{2x^4+x^3-1}{(2x^3+x^2-3x+1)(x-1)}.$$
\end{thm}

\bibliography{refs-B3geod}
\bibliographystyle{plain}

\end{document}